\documentclass{article}
\usepackage{amsmath,amsthm,amsfonts}
\newtheorem{thm}{Theorem}[section]

\newtheorem{prop}[thm]{Proposition}
\newtheorem{cor}[thm]{Corollary}
\newcommand{\prlabel}[1]{\label{#1}}
\newcommand{\prbibitem}[2]{\bibitem[#1]{#2}}
\newcommand{\R}{{\mathbb{R}}}
\newcommand{\B}{{\mathbb{P}_+T^*}}
\newcommand{\Q}{{\mathbb{Q}}}
\newcommand{\T}{{\mathbb{T}}}
\newcommand{\Z}{{\mathbb{Z}}}
\newcommand{\N}{{\mathbb{N}}}
\newcommand{\la}{{\lambda}}
\newcommand{\La}{{\mathcal{L}}}
\newcommand{\Fa}{{\mathcal{F}}}
\newcommand{\Ca}{{\mathcal{C}}}
\newcommand{\Sh}{{\text{Shape}}}

\begin{document}

\title{An obstruction to conservation of volume in contact dynamics
\thanks{Supported by the United States -Israel Binational
Science Foundation grant 1999086}
}

\author{Leonid Polterovich\\School of Mathematical Sciences\\Tel Aviv
  University\\69978 Tel Aviv \\Israel}

\date{September 26, 2000}

\maketitle

\centerline{preliminary version}
\section{Introduction and results}

A diffeomorphism of a smooth orientable
closed manifold is called {\it conservative}
if it has an invariant measure
given by a strictly positive continuous
density with respect to a smooth volume form. 
For diffeomorphisms,
there are no global topological obstructions to conservativity.
In fact a result of Moser \cite{Moser} implies
that every diffeomorphism can be isotoped to
a conservative one.
In the present note we
observe that, in contrast to this,
in the contact category such obstructions do exist.

For a closed manifold $X$ consider the space of cooriented
contact elements $\B X$ endowed with the standard contact structure
(see e.g. \cite{basic}). 
Given a contact diffeomorphism $f$
of $\B X$, write $I_f$ for the inverse of the induced 
automorphism of the
first cohomology group $H^1(\B X,\R)$
(so the mapping $f \to I_f$ is a homomorphism of the corresponding
groups).
We start with the following
result, which deals with the case when $X$ is the $n$-dimensional
torus $\T^n$.

\begin{thm}
\prlabel{torus1}
Let $f$ be a conservative contactomorphism of $\B \T^n$, $n \geq 3$. 
Then $I_f$
is a periodic automorphism: $I_f^m = \text{id}$
for some $m \in \N$.
\end{thm}

One can easily show that this result is sharp in the following sense:
every periodic automorphism of $H^1(\B \T^n,\R)$ preserving
integer cohomology can be represented as $I_f$ for some
conservative contactomorphism $f$ (see \S 2).
The case of the 2-torus $\T^2$ is slightly more involved
and will be considered separately in Theorem \ref{torus3}.
Below we give a quantitive version of Theorem \ref{torus1}
(see Theorem \ref{torus2})
as well as extension to general
manifolds $X$ (see Theorem \ref{general}) . 

The {\it lack of conservativity} (which
some times we call {\it dissipation})
can be measured
in the following way 
(cf. \cite{HK}, \S 5.1
).
Let $f$ be a contactomorphism of a closed contact manifold
$(M,\xi)$, where the contact structure $\xi$ is assumed to be 
cooriented. Pick a contact form $\lambda$ on $M$, which
means that $\xi = {\rm Kernel}\; \la$ and $\la$
agrees
with the coorientation. Note that $f^*\lambda/\lambda$
is a non-vanishing function on $M$. Consider a sequence
of real numbers
$$r_k(f,\lambda) = \max_{x \in M} {\big |} \log 
|(f^*)^{-k}\lambda/\lambda|{\big |},\; k \in \N.$$
Since any other contact form $\la'$ on $(M,\xi)$ satisfies
$\la' = F\la$ for some positive function $F$ on $M$, we have
\begin{equation}
\prlabel{bound}
|r_k(f,\la) - r_k(f,\la')| \leq \text{const}
\end{equation}
for all $k$.
A contactomorphism $f$ of $(M,\xi)$ is called {\it elliptic},
if the sequence $r_k(f,\la)$ is bounded
and {\it hyperbolic}
if there exists $c >0$ such that
\begin{equation}
\prlabel{growth}
r_k(f,\la) > ck
\end{equation}
for all sufficiently large $k \in \N$. There are
contactomorphisms such that the sequence $r_k$ 
has intermediate  growth (for instance,
$r_k \sim \log k$), but we will not focus on them
below.
When $f$ is hyperbolic, put
$\chi (f) = \sup c$, where $c$ is taken from \eqref{growth}.
In view of \eqref{bound}
the type of $f$ and the value of $\chi (f)$ for hyperbolic $f$
does not depend on the specific choice of a contact form $\la$.

\begin{prop}
\prlabel{ellip}
A conservative contactomorphism 
is elliptic.
\end{prop}
\noindent

\begin{proof}
Take a contact form $\la$ on $M$.
Assume that $f$ is a coorientation-preserving
conservative contactomorphism (the coorientation-reversing
case can be treated similary), so
$f$ preserves a volume form $F\la (d\la)^n$ for
some continuous positive function $F$. A simple calculation
shows that then $f$ preserves the continuous contact
form $F^{\frac{1}{n+1}}\la$, and it is easy to conclude that
$f$ is elliptic. 
\end{proof}

For a linear automorphism $I$ of a finite-dimensional real
vector space put
$$s(I) = \max {\big |} \log |c|{\big |},$$
where $c$ runs over all (complex) eigenvalues of $I$.
If $s(I) > 0$, the automorphism $I$ is called {\it hyperbolic}.

\begin{thm}
\prlabel{torus2}
Let $f$ be a contactomorphism of $\B \T^n$, $n \geq 2$.
Assume that $I_f$ is hyperbolic. Then $f$ is hyperbolic
and $\chi (f) \geq s(I_f)$.
\end{thm}

Let us turn to the case when $X$ is a general closed
manifold. For a diffeomorphism $e$ of $X$ 
write $\Gamma_{e}$ for the connected
component of the group of all contactomorphisms of
$\B X$ containing the canonical lift
of $e$ to $\B X$.
We are able to prove a version of Theorems
\ref{torus1} and \ref{torus2} under an additional assumption
that the contactomorphism 
$f$ lies in $\Gamma_e$. 
As a compensation, we estimate $\chi (f)$ in terms of
its action on the fundamental group rather than homology.
Denote by $\Ca$ the set of conjugacy classes of $\pi_1(X)$.
Fix a system of generators of $\pi_1(X)$. Given $\gamma \in \Ca$,
write $\ell (\gamma)$ for the minimal word length of an element
representing $\gamma$. Let $J : \Ca \to \Ca$ be the map induced by $e$.
Put 
$${\bar s} (J) = \sup _{\gamma \in \Ca} \limsup _{n \to \infty}
\frac{1}{n}\log \ell (J^n\gamma).$$
It is easy to see that this quantity does not depend on the
choice of a system of generators. We call $J$ {\it hyperbolic}
if ${\bar s} (J) > 0$.

\begin{thm}
\prlabel{general}
Assume that $J$ is hyperbolic. Then every contactomorphism
$f \in \Gamma_{e}$ is hyperbolic and $\chi (f) \geq {\bar s}(J)$.
In particular, $f$ is non-conservative.
\end{thm}

Denote by $I$ the inverse of the automorphism
of $H^1(X,\R)$ induced by $e$. It is easy to see that
${\bar s}(J) \geq s(I)$. 

\begin{cor}
\prlabel{genc}
Assume that $I$ is hyperbolic. Then every contactomorphism
$f \in \Gamma_{e}$ is hyperbolic and $\chi (f) \geq s(I)$.
\end{cor}

Here is an example where
Theorem \ref{general}  gives a better estimate than Corollary \ref{genc}.
Let $X$ be a closed oriented surface
of genus $\geq 2$. There exists a pseudo-Anosov diffeomorphism
of $X$ acting {\it trivially} on $H^1(X,\R)$ (see e.g. \cite{LP}).
On the other hand every pseudo-Anosov diffeomorphism induces
a hyperbolic map $J: \Ca \to \Ca$ (see \cite{FLP}), 
so ${\bar s}(J) > s(I) = 0$.

\medskip
\noindent
There are several possible viewpoints on the
dissipation phenomenon described above. One of them
is closely related to {\it diffusion} for symplectomorphisms
of cotangent bundles. We use it in the next section in order
to prove Theorem \ref{torus1} (and its version for the 2-torus)
and \ref{torus2}. Another approach based on the isometric
action of the contact mapping class group on the moduli
space of contact forms is presented in \S 3 where we prove
Corollary \ref{genc}. This action is studied with the
use of symplectic shapes introduced by Sikorav and Eliashberg.
In \S 4 we prove Theorem \ref{general} using 
Floer-Hofer symplectic homology
theory. 
 Finally,
\S 5 contains some speculations and open problems, and in particular
we discuss some duality between symplectic homology and shapes. 
Let us emphasize
that our decision to present a separate proof 
of Corollary \ref{genc} is related
to the fact that it is much more transparent than the
more general argument given in \S 4. 

\medskip
\noindent
{\bf Acknowledgment.} I thank Paul Biran, Dima Burago,
Yasha Eliashberg, Misha Farber and Karl Friedrich Siburg
 for very useful
discussions. I am grateful to Kai Cieliebak and Jean-Claude Sikorav for
valuable comments on the first draft of the paper. 
The isometric action of the contact mapping
class group presented in \S 4 
appeared in a stimulating  conversation with Eliashberg.
Results of this paper were presented first in KIAS (Seoul)
in Spring, 2000. I thank Yong-Geun Oh for his warm hospitality
during this visit.

\section{Symplectic diffusion}

Our proof of the theorems stated above
starts with the following {\it symplectization} procedure 
(see e.g. \cite{basic}).
For a closed manifold $X$ consider its cotangent bundle
$\pi: T^*X \to X$ endowed with the standard symplectic 
form $\Omega = dp \wedge dq$. Write $Z$ for the zero section
$\{p=0\}$, and put $T^*_0 X = T^*X \setminus Z$.
The group $\R_+$ acts naturally on $T^*_0X$ by 
$(p,q) \to (cp,q), \; c\in \R_+$, and the corresponding space
of orbits is naturally identified with $\B X$.
Write $\tau : T^*_0X \to \B X$ for the natural projection.
Let us look more attentively at this fibration.
A point $(p,q) \in T^*_0 X$ represents a contact element
\footnote{That is a tangent hyperplane to $X$}
 $ l =\text{Kernel}(p)$ at
a point $q \in X$ cooriented by $p$, so $\tau (p,q) = (l,q)$.
Further, $p$
is considered as a covector, this time, on $\B X$ at a point
$\tau(p,q)$ which sends any tangent vector $({\dot l},{\dot q})$
to $p({\dot q})$. Its cooriented kernel is precisely the
cooriented contact hyperplane on $\B X$. 
In view of this discussion,
every coorientation preserving contactomorphism $f$
of $\B X$ lifts to a unique  $\R_+$-equivariant
diffeomorphism $\phi$ of $T^*_0X$ given by
$$\phi (p,q) = ((f^*)^{-1}p,f(\tau(p,q)),$$
which in addition
turns out to be symplectic with respect to $\Omega$.
Vice versa, every $\R_+$-equivariant symplectomorphism
of $T^*_0X$ is the lift of a unique contact transformation.

Further, contact forms $\lambda$ are in a one-to-one
correspondence with smooth
sections of the bundle $T^*_0X \to \B X$.
Fix such a $\lambda$, and
for a point $z = (p,q) \in T^*_0X$ set $|z| = p/\lambda$,
where $\lambda$ is taken at the point $\tau (z)$.

The construction above allows us to translate
contact dissipation described in \S 1 into symplectic
language. The next proposition is an immediate
consequence of definitions, and its proof is left 
to the reader.

\begin{prop}
\prlabel{diff}
Let $f$ be a coorientation preserving contactomorphism
of $\B X$, and let $\phi$ be its lift to an $\R_+$-equivariant
symplectomorphism of $T^*_0X$.
\begin{itemize}
\item  $f$ is elliptic if and only there exist
constants $C > c > 0$ such that
$c|z| < |\phi^{k}(z)| < C|z|$
for every $z \in T^*_0X$ and $k \in \N$.

\item $f$ is hyperbolic if and only if there exists 
$c > 0$ 
such that 
$$\sup_{z \in T^*_0X}{\big |}\log\frac{|\phi^{k} (z)|}{|z|}{\big |} 
> ck$$ 
for all sufficiently large $k \in \N$.
In this case $\chi (f) \geq c$.
\end{itemize}
\end{prop}

\medskip
\noindent
Thus the fact that a contactomorphism
$f$ is non-elliptic can be interpreted
as {\it diffusive} behavior of its symplectic lift $\phi$,
while hyperbolicity of $f$ corresponds to exponentially fast
diffusion.

Our approach to the study of symplectic diffusion
described in \S 2 and \S 3
goes through symplectic topology of Lagrangian submanifolds
in cotangent bundles. 
Denote by $\La$ the space of all embedded
Lagrangian submanifolds $L \subset T^*X$ diffeomorphic
to $X$ such that
the projection $\pi{\big |}_L : L \to X$ induces
an isomorphism $i_L: H^1(X,\R) \to H^1(L,\R)$. 
For a Lagrangian submanifold $L \in \La$ denote by
$a(L) \in H^1(X,\R)$ the preimage of its Liouville class
\footnote{The Liouville class of a Lagrangian
submanifold $L$ is represented by the restriction of the
canonical 1-form $pdq$ to $L$, see \cite{basic}.}
under $i_L$. The role played by Lagrangian submanifolds
becomes especially transparent in the case when $X$ is the
torus $\T^n$.

\medskip
\noindent
{\bf Proof of Theorems \ref{torus1} and \ref{torus2} for $n \geq 3$:}
Assume without loss of generality that our contactomorphism
$f$ preserves the coorientation of the standard contact structure
on $\B X$ (otherwise consider $f^2$ instead of $f$). 
Since $n \geq 3$, there is a
natural identification 
$$H^1(T^*_0\T^n,\R) = H^1(\B \T^n, \R) = H^1 (\T^n,\R) = \R^n,$$
where $\R^n$ stands for a fiber of $T^*\T^n$. With this identification
$I_f$ is an automorphism of $\R^n$ which preserves
the lattice $H^1(\T^n,\Z) = \Z^n$. 
Denote by $||w||$ the Euclidean norm of a vector in $\R^n$. 
Choose a contact
form $\lambda$ on $\B \T^n$ 
such that $|(p,q)| = ||p||$ for all $(p,q) \in T^*_0\T^n$.

Let $\phi$ be the symplectic lift of $f$ to a symplectomorphism
of $T^*_0\T^n$.
For a vector 
$v \in \R^n \setminus\{0\}$ consider the "flat"
Lagrangian torus $L_v = \{p = v\}$.
Observe that
$a(\phi^{k}(L_v)) = I_f^k (v)$ for all
$ k \in \N$.

Now we use a fundamental result due to 
Gromov \cite{Gr}: given $w \in \R^n$,
every Lagrangian torus $L \in \La$ with $a(L) = w$ intersects 
$L_w$.  
Pick a sequence of vectors $v_k$ with $||v_k|| = 1$ and set
$w_k = I_f^k (v_k)$. We conclude that 
$\phi^{k} (L_{v_k})$ intersects $L_{w_k}$ for all $k \in \N$.
 
Suppose now that $I_f$ is not periodic. Since
the operator norm of $I_f^k$ goes to infinity 
\footnote{
Indeed,  the automorphism $I_f$ preserves the lattice
$H^1(X,\Z) \subset H^1(X,\R)$. Given $C > 0$, there
is only a finite number of such automorphisms with 
operator norm $\leq C$ (since
there is only a finite number of frames
in $H^1(X,\Z)$ consisting of vectors whose length does
not exceed  a given constant). Thus if the sequence
$||I_f^k||$
contains an infinite bounded subsequence then
$I_f^{k_1} = I_f^{k_2} $ for some $k_1 \neq k_2$,
so $I_f$ is periodic.} with $k$,
there exists
a sequence $v_k$ as above such that $w_k \to \infty$
as $k \to \infty$. Therefore there exists a sequence
$z_k $ with $|z_k| = 1$ such that $|\phi^{k}(z_k)| \to \infty$.
Hence $f$ is non-elliptic in view of the 
first part of Proposition \ref{diff},
and therefore $f$ is non-conservative in view of 
Proposition \ref{ellip}.
This completes the proof of Theorem \ref{torus1}.

Assume now that $I_f$ is hyperbolic. Then the sequence $v_k$ 
can be chosen in such 
a way that 
${\big |}\log ||w_k||{\big |} \geq s(I_f)k$ for all 
sufficiently large $k \in \N$.
Hence  we get a sequence of points
$z_k$ with $|z_k| = 1$ such that 
$${\big |}\log |\phi^{k} (z_k)|{\big |} \geq s(I_f)k $$
for all sufficiently large $ k \in \N$.
Now the second part of Proposition \ref{diff} yields
the hyperbolicity
of $f$, as well as the inequality $\chi (f) \geq s(f)$.
This completes the proof of Theorem \ref{torus2} for $n \geq 3$.
\qed 

\medskip
\noindent
It remains to handle the case $n=2$. 
Identify $\B \T^2$ with the 3-torus
$$\T^3 = \{p_1^2 + p_2^2 = 1\} \subset T^*\T^2.$$
Endow $\T^3$ with coordinates $(\theta, q_1,q_2)$,
where $q_1,q_2$ are coordinates in the base $\T^2$,
while $\theta$ is the angular coordinate in the fiber.
The standard contact structure on $\T^3$ is given by
$\text{Kernel} \; \la$, where
$$\la = \cos 2\pi\theta dq_1 + \sin 2\pi\theta dq_2.$$
Identify $H^1(\T^3,\R)$ with the space $\R^3$ generated
by the classes $[d\theta],[dq_1],[dq_2]$, and consider
the line $V$ generated by $[d\theta]$. It was shown
in \cite{EP0} that an automorphism $I \in SL(3,\Z)$
of $\R^3$ can be represented by a contactomorphism
(that is $I=I_f$ for some contactomorphism $f$)
if and only if $I(V) = V$. Such an automorphism acts 
in the following way:
$I([dx]) = \alpha[dx]+\beta[dy] + l[d\theta]$ and $I([dy]) = \gamma [dx]
+\delta [dy] + m[d\theta]$ where the matrix
$$ A_I = \begin{pmatrix} \alpha & \gamma \\ \beta & 
\delta \end{pmatrix}$$
belongs to $\text{GL}(2,\Z)$ and $l,m$ are arbitrary integers.
Here is an analogue of Theorem \ref{torus1} for the 2-torus.

\begin{thm}
\prlabel{torus3}
If an automorphism $I \in SL(3,\Z)$ represents a 
conservative contactomorphism of $\B \T^2$,   
the matrix $A_I$ is periodic.
\end{thm}

\medskip
\noindent
{\bf Proof of Theorems \ref{torus3} and \ref{torus2} for $n=2$:}
One should proceed exactly as in  the case $n \geq 3$ presented above
replacing $I_f$ by $A_I$ and taking into account that
$s(I)$ = $s(A_I)$ and that the Liouville class
of a flat torus $L_v$ transforms as follows:
$a(\phi (L_v)) = A_I(v)$.
\qed   

\medskip
\noindent
We conclude this section with a remark about sharpness
of our results for tori. Consider first the case
 $n \geq 3$. Every periodic automorphism $I \in \text{GL}(n,\Z)$
of $H^1(\B \T^n,\R) = \R^n$ is 
represented by a periodic contactomorphism $f$ which preserves
the coorientation. Indeed,
take the matrix $(I^T)^{-1}$ 
(where $I^T$ stands for the transposed matrix)
and consider it as
an algebraic automorphism of $\T^n = \R^n/\Z^n$. 
Define now $f$ as the canonical lift
of $(I^T)^{-1}$ to $\B \T^n$. Every periodic 
coorientation preserving contactomorphism
has an invariant contact form (use the obvious averaging procedure), and
therefore has a smooth invariant measure.

In the case $n=2$ periodicity of  the matrix $A_I$ does not imply
periodicity of the automorphism $I$. It turns out that in this case
there exist volume-preserving contactomorphisms $f$ with non-periodic
$I_f$. 

\begin{prop}
\prlabel{example}
Every automorphism $I \in \text{SL}(3,\Z)$ of $H^1(\T^3,\R)$
with periodic $A_I$ can be represented by a
volume preserving contactomorphism. 
\end{prop}

\begin{proof}
Indeed, each such $I$ 
represents the action in cohomology of a diffeomorphism
of the form $A^mB^kC$,
where $C$ is the canonical lift to $\T^3$ of a periodic 
automorphism of $\T^2$, 
$$A(\theta,q_1,q_2) = (\theta,q_1 + \theta, q_2)\;\;
\text{and} \;\;
B(\theta,q_1,q_2) = (\theta,q_1,q_2 + \theta).$$
Obviously, $C$ is homotopic to a volume preserving contactomorphism
of $\T^3$. It suffices to show that the same is valid for  $A$ and $B$.
We explain this for $A$ (the argument for $B$ is analogous).
In fact, the required volume preserving contactomorphism $f$
homotopic to $A$
is given by the following explicit formula:
$$f(\theta,q_1,q_2) =
(\theta, q_1 + \theta -\frac{1}{4\pi}\sin 4\pi\theta,
q_2 + \frac{1}{4\pi}\cos 4\pi\theta).$$
This completes the proof.
\end{proof}

\section{Contact mapping class group and shapes}

The argument presented in \S 2 does not work for general
manifolds $X$ where, for instance, it can happen that every Lagrangian
submanifold $L \in \La$ intersects the zero section $Z \subset T^*X$,
so the action of the symplectic lift of a contactomorphism
on Lagrangians
is not defined! We go round this difficulty by using
{\it symplectic shapes} introduced in \cite{Si},\cite{Si1} and 
\cite{E}.
Denote 
by $\La _0$ the subset of $\La$ consisting of all those $L$
which are Lagrangian isotopic to the zero section.  
Given a domain $U \subset T^*X$ put
$$\Sh (U) = \{a(L) \; {\big |}\; L \subset U,\; L\in \La_0\} \subset
H^1(X,\R).$$ 
The following elementary properties of shapes are important
for our purposes.
First of all, shapes are monotone in the following sense:
if $U \subset U'$ then $\Sh (U) \subset \Sh (U')$. 
Further, for $U \subset
T^*X$ and $c > 0$ one has 
$\Sh (cU) = c\Sh (U)$. 
Finally,
shapes behave nicely under symplectomorphisms $\psi$ of $T^*X$
preserving the zero section. Identify $H^1(T^*X,\R)$
with $H^1(X,\R)$ and write $I_{\psi}$ for the inverse of
the automorphism
of $H^1(X,\R)$ induced by $\psi$. Then 
$\Sh (\psi U) = I_{\psi} \Sh(U)$.
This follows from two obvious facts: $a(\psi(L)) = I_{\psi}a(L)$
for every $L \in \La_0$, and $\psi(\La_0) = \La_0$.
Here one should use that $\psi$
preserves the zero section.

Let us formulate the last, crucial, property of shapes
established by Sikorav with the use
of Lagrangian intersections theory (see \cite{Si1}, Proposition 2.6) .
Fix  a contact form $\lambda$ on $\B X$. 
Consider the subset 
\begin{equation}
\prlabel{U}
U_{\la} = \{
(p,q) \in T^*_0X \;{\big |}\; p/\la< 1\} \cup Z \subset T^*X.
\end{equation}
Obviously, $\Sh (U_{\la})$ is an open starshaped subset 
of $H^1(X,\R)$. Sikorav proved that in addition this subset is
bounded.
  
Let $(M,\xi)$ be a closed contact manifold
with a cooriented contact structure.
Denote by $\Gamma$ the group of 
all contactomorphisms of $(M,\xi)$ which preserve the 
coorientation, and by $\Gamma_0$ the connected
component of the identity.
Consider the {\it contact mapping class group}
$$G  = \Gamma /\Gamma_0 .$$

Let $\Fa$ be the space of all contact forms on $M$.
Put $W = \Fa/\Gamma_0$, and consider the action
of the contact mapping class group $G$ on $W$
given by the standard action of diffeomorphisms
on forms (cf. section 3.4 in \cite{Si1}).
It turns out that $G$ acts by isometries of a natural
pseudo-distance on $W$, and moreover this action
carries non-trivial information on contact dissipation.

The pseudo-distance, say $d$, is defined as follows. Take
two elements $a,b \in W$. Define
$$d(a,b) = \inf_{\alpha ,\beta} \max_{x \in M}
{\big |} \log \frac{\alpha}{\beta} (x) {\big |},$$
where the infimum is taken over all pairs of
contact forms
$\alpha,\beta \in \Fa$ representing the classes $a$ and $b$
respectively. 
\footnote{It would be interesting to decide whether
$d$ is a non-degenerate distance on $W$ or not.
For instance, as Kai Cieliebak pointed out,
$d$ is non-degenerate when $M$ is the circle $S^1$
equipped with
the $0$-dimensional cooriented contact structure.
In this case contact forms are simply the volume forms
on the circle while the group $\Gamma_0$ coincides
with $\text{Diff}_0(S^1)$. Since  any two volume forms
are diffeomorphic if and only if their volumes
coincide, we conclude that $W = \R_+$. Further,
$(W,d)$ is isometric to the Euclidean line $\R$,
where the isometry is given by $x \to \log x$. 
}

The group $G$ acts on $W$ by isometries of $d$.
For an element $g \in G$ consider its {\it displacement}
$$\text{disp}(g)  = \lim_{k \to +\infty} k^{-1}d(a,g^k a),\; a \in W.$$
Obviously, the limit exists and does not depend on the choice of the
point $a \in W$. 

\begin{prop}
\prlabel{disp}
Let $f \in \Gamma$ be a contactomorphism
representing a class $g \in G$.
\begin{itemize}
\item If $f$ is elliptic then for every $a \in W$ the trajectory
$\{g^ka\},\; k \in \N$ is bounded;
\item If $\text{disp}(g) > 0$ then $f$ is hyperbolic
and $\chi (f) \geq \text{disp} (g)$.
\end{itemize}
\end{prop}

\noindent
This is an immediate consequence of definitions.

Our next goal is to prove Corollary \ref{genc} with the use of this
language.
Let $M = \B X$. For a contact form
$\lambda$ consider the set $U_{\la}$ defined in \eqref{U}.
Clearly, $\Sh (U_{\la}) $ does not alter when one considers
$(f^*)^{-1}\lambda$ instead of $\lambda$ for $f \in \Gamma_0$,
so for $a \in W$ we have a well defined subset of $H^1(X,\R)$
denoted by $\Sh (a)$. This subset is open,
starshaped and bounded. Denote by $V$ the space of all open, bounded, 
starshaped domains of $H^1(X,\R)$. It carries a natural metric
$\delta (A,B) = \inf \log c$, where the infimum is taken
over all numbers $c > 1$ such that $B \subset cA$ and $A \subset cB$.
Clearly, the map $\Sh: (W,d) \to (V,\delta)$ does not increase
the (pseudo)distance. 
 
Denote by $D$ the group of all linear automorphisms
 of $H^1(X,\R)$ preserving the lattice $H^1(X,\Z)$. 
Consider the usual mapping class group
$$K = \text{Diff}(X)/\text{Diff}_0(X)$$
of the manifold $X$. We have two natural homomorphisms:
\begin{itemize} 
\item $j: K \to G$, which takes a diffeomorphism of $X$
to its canonical lift to a contactomorphism of $\B X$;
\item $I: K \to D$ which takes a diffeomorphism
to the inverse of the induced automorphism in cohomology.
\end{itemize}

\begin{prop}
\prlabel{funct}
For every $u \in K$ and $a \in W$ we have
$\Sh (j(u)a) = I(u) \Sh(a)$. 
\end{prop}

For the proof,
we need the following auxiliary fact which will be useful in the next section
as well. 
Let $\phi$ be an $\R_+$-equivariant
symplectomorphism of $T^*_0X$, and let $\psi$ be a
symplectomorphism of $T^*X$
preserving the zero section
whose restriction to $T^*_0X$ is $\R_+$-
equivariant. Assume that $\phi$ and $\psi {\big |}_{T^*_0X}$
are isotopic through $\R_+$-equivariant symplectomorphisms
of $T^*_0X$. Take any contact form $\la$ on $\B X$ and
put $U = U_{\la}$.

\begin{prop}
\prlabel{ext} 
There exists $ 0.5 > \epsilon > 0$ and a symplectomorphism
$h: T^*X \to T^*X$ such that $h \equiv \phi$ outside $0.5 U$
and $h \equiv \psi$ in $\epsilon U$.
\end{prop}

\medskip
\noindent
{\bf Proof:}
Write $\phi = \psi \theta$, where $\theta$ is a $\R_+$-equivariant
symplectomorphism of $T^*_0X$. Then there exists an isotopy 
$\theta_t, \; t \in [0;1]$ with $\theta_0 = \text{id}$
and  $\theta_1 = \theta$
consisting of
$\R_+$-equivariant symplectomorphisms of
$T^*_0X$. It suffices to show that there exists a symplectomorphism
$\sigma$ of $T^*X$ which coincides with $\theta$ outside $0.5U$ and
equals identity in $\epsilon U$. Indeed, then $h = \psi \sigma$ is as
required. To construct $\sigma
$, choose $0.01 > \epsilon > 0$ such that
$$\cup_t \theta_t (T^*X \setminus  0.4 U) \subset T^*X \setminus 
4 \epsilon U.$$ Let $B: T^*X \to \R$ be a cut off function which
vanishes on $2 \epsilon U$ and equals 1 outside $3\epsilon U$.
The path $\theta_t$ is generated by a 
degree 1 homogeneous Hamiltonian
function $A(z,t)$ on $T^*_0X$. Consider the Hamiltonian isotopy
$\sigma_t, \; t \in [0;1]$ generated by Hamiltonian $B(z)A(z,t)$.
Clearly, $\sigma = \sigma_1$ is as required.
\qed

\medskip
\noindent
{\bf Proof of Proposition \ref{funct}:}
Take a contact form $\lambda$ representing some $a \in W$.
Let $e$ be a diffeomorphism of $X$ representing $u$, and
let $f$ be a contactomorphism representing $j(u)$. Denote by $\psi$
the canonical lift of 
$e$ to a symplectomorphism of $T^*X$, and by $\phi$ the
canonical lift of $f$ to an $\R_+$-equivariant 
symplectomorphism of $T^*_0X$. Proposition \ref{ext}
yields the existence of a
symplectomorphism $h$ of $T^*X$ which coincides
with $\phi$ outside $0.5U_{\la}$ and preserves the zero section.
Looking at the action of $h$ on the
boundary of $U$ we conclude that $h(U) = U_{\beta}$,
where $ \beta =(f^*)^{-1}\lambda$
is the contact form
representing $j(u)a \in W$.
Since for every Lagrangian submanifold $L \in \La_0$ one has
$a(h(L)) = Ia(L)$ we conclude that $\Sh (j(u)a) = I \Sh (a)$.
\qed

\medskip
\noindent
{\bf Proof of Corollary\ref{genc}:}
For every $a \in W$
$$d(a,j(u)^ka) \geq \delta(\Sh (a), I(u)^k \Sh (a)).$$
Since $I$ is hyperbolic, the expression at the right hand side
grows asymptotically as $s(I)k$. Therefore,
$\text{disp}(j(u)) \geq s(I) > 0$, so
 Proposition
\ref{disp} implies that every contactomorphism $f$ representing
$j(u)$ is hyperbolic with $\chi(f) \geq s(I)$.
\qed

\medskip
\noindent
It would be interesting to explore further the isometric
action
of the mapping class group $G$ on the space $W$. Let us mention
also that one can define a canonical
isometric action of $G$ on the metric space $Z$ introduced in 
\cite{EP} and prove a direct analogue of Proposition \ref{disp}
in this case.

\section{Using symplectic homology}

Let us recall very briefly the symplectic homology theory,
see \cite{CFH, Viterbo2, Ci}. For a conjugacy class $\gamma \in \Ca$
denote by $\Lambda_{\gamma}$ the space of smooth loops $S^1 \to T^*X$
whose projection to $X$ represent $\gamma$. Let $\Phi$ be the set
of all smooth Hamiltonians $F : T^*X \times S^1 \to \R$ such that
there exist a compact subset $Q \subset T^*X$, a strictly convex 
function $\kappa: \R_+ \to \R_+$ and a contact form $\lambda$ such that
$F(p,q) = \kappa(p/\la)$ for all $(p,q) \in T^*X \setminus Q$. For
$F \in \Phi$ consider the action functional 
$$A : \Lambda_{\gamma} \to \R,\;\; A(z) = \int_z pdq - Fdt.$$
Its critical points are precisely 1-periodic orbits of the Hamiltonian
vector field $\xi$ given by $i_{\xi}\Omega = -dF$. Fix $a \in \R$,
and put $E^a_{\gamma} (F)$ to be the homology of the Floer-Morse complex
of the set $\{z \in \Lambda_{\gamma},\; A(z) < a\}$ associated to the negative
gradient flow of $A$. We work with $\Z_2$-coefficients, so $E^a_{\gamma}(F)$
is a vector space over $\Z_2$. There are  the following
natural maps between these spaces:
\begin{itemize}
\item $E^a_{\gamma}(F) \to E^b_{\gamma}(F)$ for all $F \in \Phi$ and
all pairs $a,b \in \R$ with $a < b$;
\item $E^a_{\gamma}(F_1) \to E^a_{\gamma}(F_2)$ for all $a \in \R$
and all pairs of functions $F_1,F_2 \in \Phi$ with $F_1 \leq F_2$.
\end{itemize}
Using the second natural map, one defines the symplectic homology
$E^a_{\gamma}(U)$ of a bounded domain $U \subset T^*X$ as the
(unique!)
direct limit of the spaces $E^a_{\gamma}(F)$, where $F$ runs over 
the set $\{F \in \Phi, \; F{\big |}_{U \times S^1} \leq 0\}$ 
directed by the natural
partial order $\leq$ on functions.  
The natural maps above respect the direct
limit and give rise to natural maps
\begin{itemize}
\item $E^a_{\gamma}(U) \to E^b_{\gamma}(U)$ for $  a < b$;
\item $E^a_{\gamma}(V) \to E^a_{\gamma}(U)$ for $ U \subset V$.
\end{itemize}
Define $E^{\infty}_{\gamma}(U)$ as the direct limit of $E^a_{\gamma}(U)$
when $a \to +\infty$.
Let us list some properties of symplectic homology
(all of them besides 4.3 are simple consequences of definitions).

\medskip
\noindent
{\bf 4.1.} $E^a_{\gamma}(cU) = E^{ca}_{\gamma}(U)$ for all $c \in \R_+$
and $a \in \R \cup +\infty$.

\medskip
\noindent
{\bf 4.2.} Let $h: T^*X \to T^*X$ be a symplectomorphism
of $T^*X$ which is $\R_+$ - equivariant outside a compact set. 
Denote by $J: \Ca \to \Ca$ the map induced by $h$.
Then $h$ induces an isomorphism
$$E^a_{\gamma}(U) \to E^a _{J\gamma}(h(U)).$$
   
\medskip
\noindent
{\bf 4.3 (Viterbo\cite{Viterbo2},Weber\cite{Weber}).} Let
 $U = U_{\lambda}$
be the domain associated to 
a contact form $\lambda$ on $\B X$ defined by
formula \eqref{U} of \S 3. Then $E^{\infty}_{\gamma} (U) = H_*(\Lambda_{\gamma},\Z_2)$.

\medskip
\noindent
{\bf 4.4.} Let $\lambda$ be a contact form on $\B X$.
Put
$$l_{\gamma}(\lambda) = \inf_z \int_z pdq,$$
where the infimum is taken over all closed characteristics of
the hypersurface $\{p = \lambda\} \subset T^*X$. Here each closed
characteristic is equipped with the natural orientation. Then
$E^a_{\gamma} (U_{\lambda}) = 0$ for every $a < l_{\gamma}(\lambda)$.

\medskip
\noindent
For a contact form $\lambda$ on $\B X$ define $\nu_{\gamma}(\lambda)$
as the infimum of those $a \in \R$ for which the natural map
$E^a_{\gamma}(U_{\lambda}) \to E^{\infty}_{\gamma}(U_{\lambda})$ does not vanish.
The following properties of $\nu$ are important for our purposes.
 
\medskip
\noindent
{\bf 4.5.} The quantity $\nu_{\gamma}(\la)$ is finite (see 4.3 above)
and satisfies $\nu_{\gamma}(\lambda) \geq l_{\gamma}(\lambda)$ (see 4.4 above).

\medskip
\noindent
{\bf 4.6.} If $\lambda_1 \leq \lambda_2$ then $\nu_{\gamma}(\lambda_1) \leq
\nu_{\gamma}(\lambda_2)$. Indeed, $U_{\lambda_1} \subset U_{\lambda_2}$,
and the desired statement follows from the sequence of maps
$$E^a_{\gamma}(U_{\lambda_2}) \to E^a_{\gamma}(U_{\lambda_1}) \to
E^{\infty}_{\gamma}(U_{\lambda_1})=E^{\infty}_{\gamma}(U_{\lambda_2}).$$
Here all the arrows as well as their composition are natural maps considered
above.

\medskip
\noindent
{\bf 4.7.} $\nu_{\gamma}(c\lambda) = c\nu_{\gamma}(\lambda)$ for all $c \in \R_+$.

\medskip
\noindent
{\bf 4.8.}  Let $h: T^*X \to T^*X$ be a symplectomorphism
of $T^*X$ which is $\R_+$ - equivariant outside a compact set. 
Denote by $J: \Ca \to \Ca$ the map induced by $h$. Let $\lambda_1$
and $\lambda_2$ be contact forms on $\B X$ such that
$h(U_{\lambda_1}) = U_{\lambda_2}$.
Then $\nu_{J\gamma}(U_{\la_1}) = \nu_\gamma(U_{\la_2}).$ 
This is an immediate consequence of 4.2.

\medskip
\noindent
{\bf Proof of Theorem \ref{general}:}
Take a contactomorphism $f \in \Gamma_e$ and consider
its lift $\phi$ to $T^*_0X$. Take any contact form $\lambda$
on $\B X$. Modifying $\phi$ in a neighbourhood of the zero section
and extending it to the zero section as in Proposition \ref{ext}
we get a symplectomorphism $h$ which preserves the zero section
such that $h(U_{\lambda}) = U_{(f^*)^{-1}\la}.$ Using properties
4.6-4.8 we get that 
$$d([\la],[(f^*)^{-1}\la]) \geq  {\Big |}\log \frac{\nu_{J\gamma}(\la)}{\nu_{\gamma}(\la)}
{\Big |}\;,$$
where $[\la]$ and 
$[(f^*)^{-1}\la]$ are the classes of the corresponding contact forms in the moduli
space $W$ and $d$ is the pseudo-metric on $W$ defined in \S 3.

Thus Proposition \ref{disp} combined with  4.5 yield
the inequality
$$\chi(f) \geq \sup_{\gamma} \limsup_{n \to \infty} 
\frac{1}{n}\log \nu_{J^n\gamma}(\la) \geq
\sup_{\gamma} \limsup_{n \to \infty} 
\frac{1}{n}\log l_{J^n\gamma}(\la)$$
which is valid for every contact form $\la$ on $\B X$.
Denote the expression on the right hand side by $K_{\la}$.
Take $\la$ to be the contact form associated to a Riemannian metric
on $X$. Then $l_{\gamma}(\la)$ is the length of the shortest loop
in the free homotopy class $\gamma$. It follows from \cite{Grbook}, 3.22 that
for every system of generators in $\pi_1(X)$ 
there exists $c > 1$ such that 
$$c^{-1}\ell (\gamma) \leq l_{\gamma}(\la) \leq c \;\ell (\gamma),$$
where $\ell (\gamma)$ is the minimal word length of an element
of $\pi_1(X)$ representing $\gamma$ (see \S 1).
Thus $K_{\la} = {\bar s}(J).$
This completes the proof.
\qed

\section{Discussion and open problems}

We start our discussion with the following well
known fact. For a diffeomorphism $f$ of
a manifold $M$ define its leading Lyapunov exponent
by 
$$\text{Lyap}(f) = \limsup_{k \to \infty} \frac{1}{k} \max_{x \in M}
 \big{|}\log |D_x(f^k)|\big{|},$$
where $|D_xf|$ denotes the operator norm of the differential
of $f$ at a point $x$ with respect to some fixed Riemannian metric 
on $M$.

\medskip
\noindent
{\bf Fact 5.1. }{\it A diffeomorphism
inducing a hyperbolic automorphism in 
the $m$-th (co)homology of a manifold has a strictly positive
leading
Lyapunov exponent.}

\medskip
\noindent
Our results described in \S 1 can be considered 
from two slightly different viewpoints. First, they
refine the  statement of 5.1 for contactomorphisms
of $ M = \B X$ in the case $m=1$. Indeed, the invariant $\chi(f)$
can be considered as the Lyapunov exponent
in the direction transversal to the contact structure, so
a priory it can be smaller than the leading Lyapunov exponent $\text{Lyap}(f)$,
and hence might vanish. We show that this is not the case.

From another viewpoint, Corollary \ref{genc}
is a {\it contactization} of 5.1 in the case
$m=1$.
Indeed, let $e$ be a diffeomorphism
of the base $X$ acting hyperbolically in (co)homology.
Let
$\tilde e$ be its canonical lift to the contactomorphism of $\B X$.
Then obviously $\chi({\tilde e}) = \text{Lyap}(e) > 0$. Corollary \ref{genc}
extends this inequality to the whole connected component
of $\tilde e$ in the group of contactomorphisms.

The proof of 5.1 is based on a transparent
(classical) geometric idea.
Since $f$  acts hyperbolically in homology,
the volume of the images of certain  $m$-cycles in $M$
under the
iterations of $f$ grows exponentially fast, thus
the norm of $D(f^k)$ must grow exponentially fast, so
$\text{Lyap} (f) > 0$. The approach we used in \S 4
can be considered as
a combination of this simple geometric
idea (in the case $m=1$) with Floer-Hofer symplectic homology theory
which enabled us to use symplectic
action on $T^*X$ as a substitute for the length of loops on $X$. 

It sounds likely that Theorem \ref{general} remains
true without the assumption
on the isotopy class of a contactomorphism.
There is hope that this can be done
along the lines of \S 4 but with the use of 
contact homology theory
 (see \cite{ECM,EGH}) instead of symplectic homology.
Moreover, it would be interesting to apply contact homology
to the study of dissipation on
contact manifolds which are more sophisticated
than spaces of cooriented contact elements.
For instance, one can expect that the hyperbolicity of
a contactomorphism can be detected in terms of its
action on contact homology. 
 
The link between shapes
and symplectic homology
is not totally clear at the moment, however there  is an indication that 
these notions can be considered as dual ones.
From a naive viewpoint the duality looks as follows:
shapes are responsible 
for Lagrangian submanifolds of $T^*X$, and thus are
related to closed 1-forms on $X$, while symplectic
homology deals with closed orbits on $\B X$ whose projections
define 1-cycles on $X$. An elaboration of this idea
leads to the following problem. Given a contact form $\la$
on $\B X$, symplectic homology gives rise to a function
$N: H_1(X,\Z) \to [0;+\infty)$ defined by
$N(\beta) = \inf_{\gamma} \nu_{\gamma}(\la)$, where
$\gamma$ runs over all free homotopy classes of loops
representing $\beta$. It sounds likely that one can stabilize this
function by considering
$$N_{\infty} (\beta) = \lim_{n \to \infty} \frac{1}{n}N(n\beta).$$
Extend $N_{\infty}$   to $H_1(X,\Q)$ by linearity
and define a subset $B_{\la} \subset H_1(X,\R)$ as the closure
of $\{N_{\infty} < 1\}$. Is there any relation  between $B_{\la}$ 
and $\Sh (U_{\la})$? For instance, when $\la$ is a contact form associated
to a Riemannian metric on $X$, there is a strong evidence that both
sets are convex and dual one another in the usual sense (see \cite{PS}).
Let me present
one result in this direction which holds true for every
contact form $\la$ on $\B X$. Namely,
\begin{equation}
\prlabel{dual}
\nu_{\gamma}(\la) \geq (b, \gamma)
\end{equation}
for all $\gamma \in \Ca$ and $b \in \Sh (U_{\la})$.
This translates into the language developed above
as $(b,\beta) \leq 1$ for every $b \in \Sh (U_{\la})$
and $\beta \in B_{\la}$. For the proof of inequality
\eqref{dual} assume that $L \subset U =U_{\la}$ is a
Lagrangian submanifold which is Lagrangian isotopic to the zero
section and whose Liouville class equals $b$. 
There exist $\epsilon > 0$ and a symplectomorphism $h$ of $T^*X$
such that $h(\epsilon U) \subset U$ and $[h^*(pdq)-pdq] = b$.
Moreover one can choose such an $h$ so that outside a compact
subset of $T^*X$ one has $h(p,q) = (p+ \theta(q),q)$,
where $\theta$ is a closed 1-form on $X$ representing $b$.
One can show that  $h$ induces a shift of real grading on symplectic homology:
$E^a_{\gamma}(h(V)) = E^{a -(b,\gamma)}(V)$ for all $a$ and $V$. Since $h$ 
takes $\epsilon U$
inside $U$ 
we have
$$\nu_{\gamma}(\epsilon \lambda) + (b,\gamma) \leq \nu_{\gamma}(\la).$$
Using that $\nu_{\gamma}(\epsilon \la) = \epsilon \nu_{\gamma}(\la)$
and taking $\epsilon > 0$ to be arbitrarily small we get
inequality \eqref{dual}. 

\medskip
\noindent
We conclude the paper with the following problem 
which is absolutely open. Is it possible to contactize
5.1 for $m \geq 2$? More precisely,
we ask the following

\medskip
\noindent
{\bf Question 5.2.} {\it 
Let $e$ be a 
diffeomorphism of a closed {\it simply connected } 
manifold $X$ which acts hyperbolically on $H^m(X,\R)$
for some $m \geq 2$ (for instance, $X = S^3 \times S^3$
and $m = 3$).
Let $f$ be a contactomorphism of $\B X$ which is isotopic
to the canonical lift of $e$ through contactomorphisms.
Is it true that $f$ is non-conservative/hyperbolic?}

\end{document}